\documentclass[a4paper,11pt]{amsart}
\usepackage{amssymb}
\usepackage{ifthen}
\usepackage{graphicx}
\usepackage{float}
\usepackage{caption}
\usepackage{subcaption}
\usepackage{cite}
\usepackage{amsfonts}
\usepackage{amscd}
\usepackage{amsxtra}
\usepackage{color}

\newtheorem{thm}{Theorem}
\newtheorem{cor}{Corollary}
\newtheorem{lem}{Lemma}
\newtheorem{rem}{Remark}

\newtheorem{conj}{Conjecture}
\newtheorem{prob}{Problem}

\theoremstyle{definition}
\newtheorem{defn}{Definition}[section]
\newtheorem{example}{Example}

\newenvironment{pf}[1][]{%
 \vskip 1mm
 \noindent
 \ifthenelse{\equal{#1}{}}%
  {{\slshape Proof. }}%
  {{\slshape #1.} }%
 }%
{\qed\bigskip}

\newcounter{alphabet}

\newenvironment{Thm}[1][]{\refstepcounter{alphabet}%
\bigskip%
\noindent%
{\bf Theorem \Alph{alphabet}}%
\ifthenelse{\equal{#1}{}}{}{ (#1)}%
{\bf .} \itshape}{\vskip 8pt}

\def\be{\begin{equation}}
\def\ee{\end{equation}}

\newcommand{\bee}{\begin{enumerate}}
\newcommand{\eee}{\end{enumerate}}

\newcommand{\blem}{\begin{lem}}
\newcommand{\elem}{\end{lem}}
\newcommand{\bthm}{\begin{thm}}
\newcommand{\ethm}{\end{thm}}
\newcommand{\bcor}{\begin{cor}}
\newcommand{\ecor}{\end{cor}}
\newcommand{\beg}{\begin{example}}
\newcommand{\eeg}{\end{example}}
\newcommand{\begs}{\begin{examples}}
\newcommand{\eegs}{\end{examples}}
\newcommand{\bdefe}{\begin{defn}}
\newcommand{\edefe}{\end{defn}}
\newcommand{\bprob}{\begin{prob}}
\newcommand{\eprob}{\end{prob}}
\newcommand{\bques}{\begin{ques}}
\newcommand{\eques}{\end{ques}}
\newcommand{\bei}{\begin{itemize}}
\newcommand{\eei}{\end{itemize}}
\newcommand{\bcon}{\begin{conj}}
\newcommand{\econ}{\end{conj}}
\newcommand{\bcons}{\begin{conjs}}
\newcommand{\econs}{\end{conjs}}
\newcommand{\bprop}{\begin{propo}}
\newcommand{\eprop}{\end{propo}}
\newcommand{\br}{\begin{rem}}
\newcommand{\er}{\end{rem}}
\newcommand{\brs}{\begin{rems}}
\newcommand{\ers}{\end{rems}}
\newcommand{\bo}{\begin{obser}}
\newcommand{\eo}{\end{obser}}
\newcommand{\bos}{\begin{obsers}}
\newcommand{\eos}{\end{obsers}}
\newcommand{\bpf}{\begin{pf}}
\newcommand{\epf}{\end{pf}}
\newcommand{\ba}{\begin{array}}
\newcommand{\ea}{\end{array}}
\newcommand{\beq}{\begin{eqnarray}}
\newcommand{\beqq}{\begin{eqnarray*}}
\newcommand{\eeq}{\end{eqnarray}}
\newcommand{\eeqq}{\end{eqnarray*}}

\newcommand{\ds}{\displaystyle}

\newcounter{minutes}
\divide\time by 60
\newcounter{hours}
\multiply\time by 60 \addtocounter{minutes}{-\time}
\begin{document}

\bibliographystyle{amsplain}
%

\title[Generalized versions of Lipschitz conditions on the modulus of holomorphic functions ]{Generalized versions of Lipschitz conditions on the modulus of holomorphic functions}

\def\thefootnote{}
\footnotetext{ \texttt{\tiny File:~\jobname .tex,
          printed: \number\day-\number\month-\number\year,
          \thehours.\ifnum\theminutes<10{0}\fi\theminutes}
} \makeatletter\def\thefootnote{\@arabic\c@footnote}\makeatother


\author[S. Ponnusamy]{Saminathan Ponnusamy
}
\address{
S. Ponnusamy, Department of Mathematics,
Indian Institute of Technology Madras, Chennai-600 036, India.
}
\email{samy@iitm.ac.in}

\author[R. Vijayakumar]{Ramakrishnan Vijayakumar}
\address{
R. Vijayakumar, Department of Mathematics,
Indian Institute of Technology Madras, Chennai-600 036, India.
}
\email{mathesvijay8@gmail.com}

\subjclass[2010]{Primary: 30C80, 30H05, 32A10 ; Secondary: 30G30, 46B20, 46E15, 46E40
}
\keywords{Holomorphic functions, Schwarz-Pick lemma, Lipschitz classes}
\begin{abstract}
In this paper, we establish Lipschitz conditions for the norm of holomorphic mappings between the unit ball $\mathbb{B}^n$ in $\mathbb{C}^n$ and $X,$
a complex normed space. This extends the work of  Djordjevi\'{c} and Pavlovi\'{c}. 
\end{abstract}

\maketitle
\pagestyle{myheadings}
\markboth{
S. Ponnusamy, R. Vijayakumar}{Generalized versions of Lipschitz conditions on the modulus of holomorphic functions}

\section{Introduction and Preliminaries}
Denote by $\mathbb{C}^n$, the $n$-dimensional complex Hilbert space with the inner product and the norm given by
$ \langle z, w \rangle:=\sum_{j=1}^{n}z_{j}\overline{w_{j}}$  and $\|z\|:=\sqrt{\langle z, z \rangle},$
where $z, w \in \mathbb{C}^n$, respectively. Write ${\mathbb{B}}^{n} :=\{z \in {\mathbb{C}}^{n}:\, \|z\|<1\}$ for  the open unit ball in ${\mathbb{C}}^n$
so that ${\mathbb{B}}^{1} =:\mathbb{D}$ denotes the open unit disk in $\mathbb{C}$. If $V$ and $W$ are two normed spaces and $U\subset V$ is open, then the
Fr\'{e}chet derivative of a holomorphic mapping $ f: \,U \rightarrow W$  is defined to be the unique linear map
$ A= f'(z) :\, V \rightarrow W$ such that
$$
f(z+h)=f(z)+f'(z)\cdot h+ o(\|h\|^{2})
$$
for $h$ near the origin of $V$. The norm of such a map is defined by
$ \| A \| = \sup_{\|z\|=1} \|Az\|.$

In 1975,   Globevnik \cite{Glob-75} introduced the notion of uniform c-convexity and proved that $L^{1}$-space possesses this property. Namely, a
complex normed space $X$ is said to be \emph{uniformly $c$-convex} if there exists a positive increasing function
$\Omega(\delta)~(\delta > 0)$ with $\Omega(0^{+})=0$ such that for all $x,y \in X$ and $\delta > 0$ there holds the implication
\beqq
\max_{|\lambda|\leq 1\atop \|x\|=1} \|x+\lambda y\| \leq 1+\delta \implies \|y\| \leq \Omega(\delta).
\eeqq
The smallest of the functions $\Omega$ is denoted by $\Omega_{X},$ i.e.,
\beqq
\Omega_{X}(\delta):= \sup\{\|y\|:\,  ~\max_{|\lambda|\leq 1 \atop  \|x\|=1 } \|x+\lambda y\| \leq 1+\delta \}.
\eeqq
As mentioned in \cite{DjoPav-08}, it can be easily seen that
$$
\Omega_{\mathbb{C}}(\delta)=\delta ~\mbox{ and }~
\Omega_{H}(\delta)=\sqrt{\delta(2+\delta)},
$$
where $H$ is a Hilbert space of dimension at least two.	

As in \cite{Dya-97}, we call a function $\omega:\, [0, \infty ) \rightarrow \mathbb{R}$ a $\emph{majorant}$ if $\omega$ is continuous, increasing, $\omega(0)=0,$ and
$ t^{-1}\omega(t)$ is nonincreasing on $(0, \infty ).$
If, in addition, there is a constant $C(\omega) > 0$ such that
$$
\int\limits_{0}^{\delta} \frac{\omega(t)}{t} dt + \delta \int\limits_{\delta}^{\infty} \frac{\omega(t)}{t^{2}} dt \leq C(\omega) \cdot \omega(\delta)
$$	
whenever $0< \delta <1,$ then we say that $\omega$ is a $\emph{regular majorant}.$
	
Then the space $\mathbf{Lip}(\omega, G, X),$ where $G$ is bounded subset of $\mathbb{C}^n,$ is defined to be the set of those functions $g:\, G \rightarrow X$ for which
$$ \|g(z)-g(w)\| \leq c \cdot \omega(\|z-w\|),$$
where $c$ is a constant. If $\omega(t)= t^{\alpha}$ for some $ \alpha \in (0,1],$ then we write $\mathbf{Lip}(\omega, G, X)= \Lambda_{\alpha}(G, X).$
If $X$ is uniformly $c$-convex, then $\Omega_{X}$ is a majorant (cf. \cite{DGT-84}).
A majorant $\omega$ is said to be a \emph{Dini majorant} if $ \int_{0}^{1}\frac{\omega(t)}{t} \,dt < \infty.$
For a Dini majorant, we define the $majorant~ \widetilde{\omega}$ by
$$ \widetilde{\omega}(t) =\int_{0}^{t}\frac{\omega(x)}{x} \,dx =\int_{0}^{1}\frac{\omega(tx)}{x} \,dx.
$$
A majorant $\omega$ is said to be \emph{fast} \cite{Dya-04} if
\beqq
\int_{0}^{\delta} \frac{\omega(t)}{t}\, dt \leq ~\mbox{const} \cdot \omega(\delta),~ 0< \delta < \delta_{0},
\eeqq
for some $ \delta_{0} > 0.$ (Of course, if $\omega$ is fast, then it is a Dini majorant).

Dyakonov \cite{Dya-97} gave some characterizations of the holomorphic functions of class $\Lambda_{\omega}(\mathbb{D}, \mathbb{C})$ in terms of their moduli. 

\begin{Thm}{\rm (\cite{Dya-97})}  \label{Po-Vi5-ThmA}
Let $\omega $ be a regular majorant. A function $f$ holomorphic in $\mathbb{D}$ is in $\Lambda_{\omega}(\mathbb{D}, \mathbb{C})$ if and only if so is its modulus $|f|.$
\end{Thm}

The main ingredient in Dyakonov's proof is a very complicated. However, Pavlovic \cite{Pav-99} gave a simple proof of Theorem~A. 
The proof uses only the basic lemmas of \cite{Dya-97} and the Schwarz lemma, and is therefore considerably shorter than that of \cite{Dya-97}. 
However, Theorem~A 
does not extend to $\mathbb{C}^{k}$- valued functions ($k \geq 2$). So we have to consider functions with
additional properties (see Theorems \ref{Po-Vi5-thm3} and \ref{Po-Vi 5 thm4}).

In \cite{DjoPav-08}, Djordjevi\'c and Pavlovi\'c extended to vector-valued functions of a theorem of Dyakonov \cite{Dya-97} on Lipschitz
conditions for the modulus of holomorphic functions.
Therefore, it is natural for us to extend this result for holomorphic functions on $\mathbb{B}^n.$
Very recently, Kalaj \cite{Kal-18} established a Schwarz-Pick type inequality for holomorphic mappings between unit balls $\mathbb{B}^n$ and
$\mathbb{B}^m$ in the corresponding complex spaces.
	
\begin{Thm}
{\rm (\cite[Theorem 2.1]{Kal-18})} \label{Po-Vi 5 Kalaj Thm}
If $f$ is a holomorphic mapping of the unit ball
${\mathbb{B}}^{n} \subset {\mathbb{C}}^{n}$  into ${\mathbb{B}}^{m}\subset {\mathbb{C}}^{m},$ then for $z\in \mathbb{B}^n$ we have
$$
\|f'(z)\| \leq \left\{
\begin{array}{rl} \ds \frac{\sqrt{1-\|f(z)\|^2}}{1-\|z\|^2} &  ~\mbox{ for $m\geq 2$},\\[3mm]
\ds \frac{1-\|f(z)\|^2}{1-\|z\|^2}  &~\mbox{for $m = 1$}.
\end{array}
\right .
$$
\end{Thm}

In \cite{DaiPan-15}, Dai and Pan proved the following theorem which establishes a Schwarz-Pick type estimates for gradient of the modulus of
holomorphic mappings.

\begin{Thm}
{\rm (\cite[Theorem 1]{DaiPan-15})} \label{Po-Vi5-ThmC}
Let $f:\, \mathbb{B}^n \rightarrow \mathbb{B}^m$ be a holomorphic mapping. Then
\beqq
\big |\nabla \|f\|(z)\big | \leq \frac{1-\|f(z)\|^{2}}{1-\|z\|^{2}} ~\mbox{for~ $z\in \mathbb{B}^n$.}
\eeqq
\end{Thm}
For a holomorphic mapping $f:\, \mathbb{B}^n \rightarrow \mathbb{B}^m,$ we have
\beq\label{Eqn 1}
\big |\nabla \|f\|(z)\big | = \frac{1}{\|f(z)\|}\bigg\|\bigg(\bigg\langle \frac{\partial f(z)}{\partial z_{1}},f(z)  \bigg\rangle,\ldots, \bigg\langle \frac{\partial f(z)}{\partial z_{n}},f(z)  \bigg\rangle\bigg)\bigg\|~\mbox{ if $f(z) \neq 0$}.
\eeq

\section{The main results}

 \begin{thm}\label{Po-Vi5-thm1}
 Let $X$ be uniformly $c$-convex and $f:\,\mathbb{B}^n \rightarrow X$ be a holomorphic function satisfying
\beq \label{Po-Vi5 - Eqn 2}
\big |\|f(z)\|-\|f(w)\|\big| \leq c \|z-w \|^{\alpha}  ~\mbox{for $z, w\in \mathbb{B}^n$},	
\eeq
where $c \geq 0$ and $\alpha \in [0, 1]$ are constants. Then
\beq\label{Eqn 2}
\|f'(z)\| \leq 2K \dfrac{\Omega_{X}(c K^{-1}(1-\|z\|)^{\alpha})}{1-\|z\|}  ~\mbox{for $z\in \mathbb{B}^n$},	
\eeq
where $K=\|f(0)\| + c.$ Especially, if $\|f(0)\|=1,$ then 	
\beq \label{PoVi5- Eqn 3}
\|f'(z)\| \leq 2(1+c) \dfrac{\Omega_{X}(c (1-\|z\|)^{\alpha})}{1-\|z\|}  ~\mbox{for $z\in \mathbb{B}^n$}.	
\eeq
\end{thm}

\begin{thm}\label{Po-Vi5-thm2}
Let $X$ be uniformly $c$-convex such that $\Omega_{X}$ is a Dini majorant and $f:\,\mathbb{B}^n \rightarrow X$ be a holomorphic function
such that the function $\|f(z)\|$ belongs to $\Lambda_{\alpha}(\mathbb{B}^n, \mathbb{R})$ for some $\alpha \in (0, 1].$ Then
$f \in \mathbf{Lip}(\overline{\omega}_{\alpha}, \mathbb{B}^n, X),$ where
$ \overline{\omega}_{\alpha} (t)= \widetilde{\Omega}_{X}(t^{\alpha}). $

In particular, the function $f$ is uniformly continuous on $\mathbb{B}^n$ that has a continuous extension to the closed disk.
\end{thm}

\begin{cor}
If $\Omega_{X}$ is fast and $f:\, \mathbb{B}^n \rightarrow X$ is a holomorphic function such that the function $\| f(z) \|$ belongs to
$\Lambda_{\alpha}(\mathbb{B}^n, \mathbb{R})$ for some $\alpha \in (0, 1].$ Then $f \in \mathbf{Lip}(\omega_{\alpha}, \mathbb{B}^n, X),$
where $ \omega_{\alpha} (t)= \Omega_{X}(t^{\alpha}).$ 	
\end{cor}

Taking $n=1$ and $X=\mathbb{C},$ we get the following result of Dyakonov \cite{Dya-97}.

\begin{cor}
If $f:\, \mathbb{D} \rightarrow \mathbb{C}$ is a holomorphic function such that $|f|$ belongs to $\Lambda_{\alpha}(\mathbb{D}, \mathbb{R})$
for some $\alpha \in (0, 1].$ Then  $f$ belongs to $\Lambda_{\alpha}(\mathbb{D}, \mathbb{C}).$
\end{cor}

\begin{thm}\label{Po-Vi5-thm3}
Let $0 < \alpha \leq 1$ and  $f:\, \mathbb{B}^n \rightarrow \mathbb{C}^m$ be a holomorphic function such that
\beq\label{Eqn 3}
\|f'(z)\| \|f(z)\| \leq K \bigg\|\bigg(\bigg\langle \frac{\partial f(z)}{\partial z_{1}},f(z)  \bigg\rangle,\ldots, \bigg\langle \frac{\partial f(z)}{\partial z_{n}},f(z)  \bigg\rangle\bigg)\bigg\|~\mbox{ for $z \in \mathbb{B}^n$},
\eeq
where $K$ is a constant independent of $z.$  Then  $f \in \Lambda_{\alpha}(\mathbb{B}^n, \mathbb{C}^m)$ if and only if $\|f\| \in \Lambda_{\alpha}(\mathbb{B}^n, \mathbb{R}).$
\end{thm}


\begin{thm}\label{Po-Vi 5 thm4}
If $f:\, \mathbb{B}^{n} \rightarrow \mathbb{C}^{m}, m\geq2,$ is holomorphic and if $\|f\| \in \Lambda_{\alpha}(\mathbb{B}^n, \mathbb{R})$
for some $\alpha \in (0, 1],$ then we have $f \in \Lambda_{\alpha/2}(\mathbb{B}^n, \mathbb{C}^m).$
\end{thm}

The case $n=1$ of Theorems \ref{Po-Vi5-thm3} and \ref{Po-Vi 5 thm4} gives results of Pavlovi\'{c} \cite{Pav-11}.
\section{Proofs of the Theorems}

Theorem \ref{Po-Vi5-thm1} is a direct consequence of the following lemma.

\begin{lem}
If $f:\,\mathbb{B}^n \rightarrow X$ is a holomorphic function satisfying the condition
\beq\label{Eqn 4}
\big |\|f(z)\|-\|f(w)\|\big | \leq c (1- \|z\|)^{\alpha}  ~\mbox{whenever  $\|w-z\| \leq 1-\|z\|$},	
\eeq
then there holds \eqref{Eqn 2}.	
\end{lem}
\begin{pf}
Fix $z \in \mathbb{B}^n$ with $f(z) \neq 0,$ and let $L \in X^{*},$ $\|L\|=1,$ where $X^{*}$ is the dual of $X.$ Consider the scalar valued function
$$ \phi (z)=L\circ f(z),
$$
and introduce the following set for the given $z \in \mathbb{B}^n,$
$$D_{z}:= \{w \in {\mathbb{C}}^{n}:\, \|w-z\|<1-\|z\| \} ~ \mbox{and}~ M_{z}:= \sup\{\|f(w)\|:\, w\in D_{z}\}.
$$
If $z=0$ and $M_{0}=1,$ then the Schwarz-Pick lemma (see Theorem~B) 
gives
$$ |\phi'(0)| \leq 1- |\phi (0)|^{2} \leq 2 (1-|\phi(0)|),$$
which is our inequality in this special case. The general case follows by applying the special case to the function $\Phi$ defined by
$$ \Phi(\zeta)= \dfrac{\phi(z + \zeta (1- \|z\|))}{M_{z}} ~ \mbox{for}~ \zeta \in \mathbb{B}^n.$$
We have
\beqq
\frac{1}{2}(1-\|z\|)|\phi'(z)|+ |\phi(z)| \leq M_{z} ~\mbox{for $z \in \mathbb{B}^n $}.
\eeqq
Since $\phi'(z)= L(f'(z)),$ we see that
$$
(1-\|z\|) |L(f'(z)/2)|+ |L(f(z))| \leq M_{z}.
$$
Hence, for every $\lambda \in \mathbb{D},$ we obtain
$$
 |\lambda (1-\|z\|) L(f'(z)/2) + L(f(z))| \leq M_{z}.
$$
Since this holds for every $L$ of norm 1, by taking the supremum over all $L$ with $\| L \|=1$ and by applying the Hahn-Banach theorem, we get
$$
 \bigg\| \lambda \dfrac{(1-\|z\|) f'(z)}{2} + f(z)\bigg\| \leq M_{z},\ \ \mbox{i.e.,}\ \
\bigg\|\dfrac{f(z)}{\|f(z)\|} + \lambda \dfrac{(1-\|z\|) f'(z)}{2 \|f(z)\|} \bigg\| \leq \dfrac{M_{z}}{\|f(z)\|}.
$$
Now denoting
$$ x= \dfrac{f(z)}{\|f(z)\|},~ y=\dfrac{(1-\|z\|) f'(z)}{2 \|f(z)\|} ~\mbox{and}~ \delta= \dfrac{M_{z}-\|f(z)\|}{\|f(z)\|},
$$
we see from the definition of $\Omega_{X}$ that
\beq\label{Eqn 5}
(1-\|z\|)\|f'(z)\| \leq 2 \|f(z)\| \Omega_{X}\left(\dfrac{M_{z}-\|f(z)\|}{\|f(z)\|}\right).
\eeq
Therefore by \eqref{Eqn 4} and \eqref{Eqn 5}, we obtain that
\beqq
(1-\|z\|)\|f'(z)\| \leq 2 \|f(z)\| \Omega_{X}\left(\dfrac{c (1- \|z\|)^{\alpha}}{\|f(z)\|}\right).
\eeqq
Now \eqref{Eqn 2} follows from the fact that $\Omega_{X}(t)/t$ is a decreasing function and the inequality $\|f(z)\| \leq K.$
The proof is complete.
\end{pf}

\begin{lem}
If a $C^{1}$-function $u:\, \mathbb{B}^n \rightarrow \mathbb{R}$ satisfies
\beqq
\| \nabla u(z) \| \leq \frac{\omega(1-\|z\|)}{1-\|z\|}~\mbox{for~ $z\in \mathbb{B}^n,$}
\eeqq
where $\omega$ is a Dini majorant, then
$$| u(a) - u(b)| \leq 3~ \widetilde{\omega}(\| a-b \|), ~\mbox{for ~ $a, b \in \mathbb{B}^n.$}
$$
\end{lem}
\begin{pf}
We begin the proof with the following observation: $\omega \leq \widetilde{\omega}.$ In fact, we let $t_0 \in (0, \infty).$
Since $\frac{\omega(t)}{t}$ is decreasing on $(0, \infty),$ we have
$$
\frac{\omega(t_0)}{t_0} \leq  \frac{\omega(t_0 x)}{t_0 x} \ \mbox{for $x \in (0,1]$}.
$$
Integrating on both sides of the last inequality from $0$ to $1$, we obtain by definition of $\widetilde{\omega}$ that
$\omega(t_0) \leq \widetilde{\omega}(t_0)$.

Let $ \|a \| \leq \|b \| \leq 1.$ By Lagrange's mean-value theorem,
$$	| u(a) - u(b)| \leq   \|\nabla u(c)\| \| a-b \|,$$
where $c = (1- \lambda) a + \lambda b$ for some $ \lambda \in (0, 1).$ Since $\|c \| \leq \|b \|$ and $ \omega(t)/ t$ decreases, we see that
$$ \frac{\omega(1-\|c\|)}{1-\|c\|}~\leq \frac{\omega(1-\|b\|)}{1-\|b\|}
$$
and hence,
$$|u(a)- u(b)| \leq \omega (\|a-b\|) \leq  \widetilde{\omega}(\| a-b \|),
$$
under the condition $\|a-b \| \leq 1- \|b\|.$	
	
If $1- \|b\| \leq 	\|a-b \| \leq 1-\|a \|,$ then
$$| u(a)- u(b) | \leq | u(a)- u(b') |+ | u(b')- u(b) |,
$$
where $b'= \frac{(1- \delta)b}{\|b\|}$ and $\delta= \|a-b\|.$ Using the Lagrange's mean-value theorem as above we get
\beqq
|u(a)- u(b')| &\leq&  \frac{\omega(1-\|b'\|)}{1-\|b'\|} \|a-b'\|
= \frac{\omega(\delta)}{\delta} \|a-b'\|
\leq \omega(\delta) \leq \widetilde{\omega}(\delta).
\eeqq	
In the case of $|u(b')-u(b)|,$ we have
\beqq	
|u(b')-u(b)|&\leq& \int_{\|b'\|}^{\|b\|} \frac{\omega(1-t)}{1-t}\, dt
\leq \int_{1-\delta}^{1} \frac{\omega(1-t)}{1-t} \,dt =  \widetilde{\omega}(\delta).
\eeqq
Finally, if $ \delta > 1-\|a\|,$ we use the inequality 	
$$| u(a)- u(b) | \leq | u(a)- u(a') |+ | u(a')- u(b') |+ |u(b')-u(b)|,
$$
where $a'= \frac{(1- \delta)a}{\|a\|},$ and then proceed in a similar way as above, using the inequality $\|a'- b'\| \leq \|a - b\|$.
\end{pf}

Lemma \ref{Po-vi5-lem3} can easily be proved by applying the previous lemma to the functions ${\rm Re}\,(L \circ f(z))$ and ${\rm Im}\,(L \circ f(z)),$ where
$L \in X^{*}$ and $\| L \|=1.$

\begin{lem}\label{Po-vi5-lem3}
If $f$ is an $X$-valued holomorphic function in $\mathbb{B}^n$ and satisfies the condition
$$\|f'(z)\| \leq \frac{\omega(1-\|z\|)}{1-\|z\|}~\mbox{ for $z\in \mathbb{B}^n$},
$$
where $\omega$ is a Dini majorant, then $f \in \mathbf{Lip}(\widetilde{\omega}, \mathbb{B}^n, X).$
\end{lem}

\subsection{Proof of Theorem \ref{Po-Vi5-thm2}}
Let $f$ satisfy the hypotheses of the theorem. Then
$$\| f'(z)/ 2 K\| \leq 	\frac{\omega(1-\|z\|)}{1-\|z\|},
$$
by Theorem \ref{Po-Vi5-thm1}, where	
$ \omega(t) = \Omega_{X}(c K^{-1} t^{\alpha}).$ But a simple calculation shows that	
$ \widetilde{\omega}(t) = {\alpha}^{-1}\widetilde{\Omega}_{X}(c K^{-1} t^{\alpha})$ and
so we can appeal to Lemma \ref{Po-vi5-lem3} to conclude the proof.
\hfill $\Box$

\subsection{Proof of Theorem \ref{Po-Vi5-thm3}}
The `only if' part is trivial. Assume that $\|f(z)\| \in \Lambda_{\alpha}(\mathbb{B}^n, \mathbb{R})$ and we
proceed as in Theorem \ref{Po-Vi5-thm1}. Fix $z \in \mathbb{B}^n$ with $f(z) \neq 0,$ and consider the following sets for a given $z \in \mathbb{B}^n,$
$$ D_{z}:= \{w \in {\mathbb{C}}^{n}:\, \|w-z\|<1-\|z\| \} ~ \mbox{ and }~ M_{z}:= \sup\{\|f(w)\|:\, w\in D_{z}\}.
$$
If $z=0$ and $M_{0}=1,$ Theorem~C 
gives
$$ |\nabla \|f\|(0)| \leq 1- \|f(0)\|^{2} \leq 2(1-\|f(0)\|).
$$
Therefore, from \eqref{Eqn 3} and the formula \eqref{Eqn 1}, we have that
$$ \|f'(0)\| \leq 2K (1-\|f(0)\|),
$$
which is our inequality in this special case. The general case follows by applying the special case to the function $F$ defined by
\beq\label{Eqn 6}
 F(\zeta)= \dfrac{f(z + \zeta (1- \|z\|))}{M_{z}} ~ \mbox{ for }~ \zeta \in \mathbb{B}^n,
\eeq
and obtain
\beq\label{Eqn 7}
\frac{1}{2K}(1-\|z\|)\|f'(z)\|+ \|f(z)\| \leq M_{z} ~\mbox{ for $z \in \mathbb{B}^n $}.
\eeq
Since $\|f\| \in \Lambda_{\alpha}(\mathbb{B}^n, \mathbb{R}),$ we have
$$ \|f(w)\| - \|f(z)\| \leq c \|w-z\|^{\alpha} \leq c (1-\|z\|)^{\alpha},
$$
for $z \in \mathbb{B}^n$ and $w \in D_{z}.$ Taking the supremum over all $w \in D_{z}$ and then using the inequality
\eqref{Eqn 7}, we get
$$ \|f'(z)\| \leq C \frac{\omega(1-\|z\|)}{1-\|z\|},
$$
where $C$ is a constant and $\omega(t)= t^{\alpha}.$
The desired conclusion follows from Lemma \ref{Po-vi5-lem3}.
\hfill $\Box$

\subsection{Proof of Theorem \ref{Po-Vi 5 thm4}}
Let $z \in \mathbb{B}^{n}$ and proceed the steps as in the above proof.
If $z=0$ and $M_{0}=1,$ then the higher dimensional version of Schwarz-Pick lemma (Theorem~C) 
gives
$$ \|f'(0)\| \leq \sqrt{1-\|f(0)\|^2} \leq \sqrt{2} \sqrt{1-\|f(0)\|},$$ which is our inequality in this special case. The general case follows by applying the special case to the function $F$ defined by \eqref{Eqn 6}. Indeed, we obtain
\beq\label{Eqn 8}
(1-\|z\|)\|f'(z)\| \leq c \sqrt{M_{z}-\|f(z)\|},
\eeq
for some constant $c$. Since $\|f\| \in \Lambda_{\alpha}(\mathbb{B}^n, \mathbb{R}),$ we have
$$ \|f(w)\| - \|f(z)\| \leq c \|w-z\|^{\alpha} \leq c (1-\|z\|)^{\alpha},$$
for $z \in \mathbb{B}^n$ and $w \in D_{z}.$ Taking the supremum over $w \in D_{z}$ and then using the inequality
\eqref{Eqn 8}, we get
$$ \|f'(z)\| \leq C \frac{\omega(1-\|z\|)}{1-\|z\|},
$$
where $C$ is a constant and $\omega(t)= t^{\alpha/2}.$
Now the result follows from Lemma \ref{Po-vi5-lem3}.
\hfill $\Box$

\begin{rem}
The index $\alpha/2$ in Theorem \ref{Po-Vi 5 thm4} is optimal as demonstrated by the following example (see \cite{Pav-11}). Consider the  function $f:\,\mathbb{D} \rightarrow \mathbb{C}^2$ by $f(z)= (1, (1-z)^{\alpha/2}),\ 0 < \alpha \leq 1.$ We have
\beqq
\bigg|\|f(z)\|-\|f(w)\| \bigg| &=& \bigg|\sqrt{\|1-z\|^{\alpha}+1} - \sqrt{\|1-w\|^{\alpha}+1}\bigg|\\
&\leq& \bigg|\|1-w\|^{\alpha} - \|1-z\|^{\alpha}\bigg| \leq \|z-w\|^{\alpha},
\eeqq
while
$\| f(1)-f(r) \|=(1-r)^{\alpha/2}, 0<r<1.$
This shows that the index $\alpha/2$ is optimal.
\end{rem}


\section{Concluding Remarks}
As mentioned in \cite{DjoPav-08}, the inequality \eqref{PoVi5- Eqn 3} is in a sense optimal for the case $n=1$. To see this, let $\omega(t) >0$ be an arbitrary increasing
function on $(0, \infty)$ such that $\omega(0^+)=0.$ We say that a Banach space $X$ has the property $\mathcal{L}(\omega, \alpha),$ if the following holds:
{\em For every $c \in (0,1)$ and every analytic function $f:\, \mathbb{D} \rightarrow X$ with $\|f(0)\|=1$, the  inequality \eqref{Po-Vi5 - Eqn 2} implies
that
$$\|f'(\lambda)\| \leq \dfrac{\omega(c (1-|\lambda|)^{\alpha})}{1-|\lambda |} \ \  ~\mbox{for $\lambda \in \mathbb{D}$}.
$$
}
It is well-known that, if the Banach space $X$ has the property $\mathcal{L}(\omega,\alpha)$  (see \cite[Proposition 10]{DjoPav-08}), then $X$ is uniformly $c$-convex and $\Omega_{X}(\delta) \leq B \omega(\delta)$ for $0< \delta <1,$ where $B$ is a constant. This result is to emphasize the fact that $\|f(0)\|=1$ provides condition for uniformly $c$-convexity of the Banach space $X.$

\end{document}